\documentclass[twocolumn]{autart}    
\newcommand*{\QEDA}{\hfill\ensuremath{\blacksquare}}
\usepackage{cite}
\usepackage{amsmath,amssymb,amsfonts}
\usepackage{textcomp}
\usepackage{color}
\usepackage{bm}
\newtheorem{assumption}{$Assumption$}
\usepackage{algorithm}
\usepackage{algorithmic}
\usepackage{savesym}
\savesymbol{AND}
\usepackage{multirow} 
\usepackage{graphicx}          

\begin{document}

\begin{frontmatter}
	
	\title{
	A unitary distributed subgradient method for multi-agent optimization with different coupling sources
} 
	
	\thanks[footnoteinfo]{This work was supported by the Natural Sciences and Engineering Research Council of Canada (NSERC). Some preliminary results have been reported in the conference version of this paper \cite{liu2019distributed}.}
	
	\author[uvic]{Changxin Liu}\ead{chxliu@uvic.ca},    
	\author[npu]{Huiping Li}\ead{lihuiping@nwpu.edu.cn}, 
	\author[uvic]{Yang Shi}\ead{yshi@uvic.ca}
	
	\address[uvic]{Department of Mechanical Engineering, University of Victoria, Victoria, B.C., Canada, V8W 3P6}  
	\address[npu]{School of Marine Science and Technology, Northwestern Polytechnical University, Xi'an, 710072, P. R. China}


\begin{abstract}                          
Distributed optimization techniques offer high quality solutions to various engineering problems, such as resource allocation and distributed estimation and control.
In this work, we first consider distributed convex constrained  optimization problems where the objective function is encoded by multiple local and possibly nonsmooth objectives privately held by a group of agents, and propose a distributed subgradient method with double averaging (abbreviated as ${\rm DSA_2}$) that only requires peer-to-peer communication and local computation to solve the global problem. The algorithmic framework builds on dual methods and dynamic average consensus; the sequence of test points is formed by iteratively minimizing a local dual model of the overall objective where the coefficients, i.e., approximated subgradients of the objective, are supplied by the dynamic average consensus scheme. We theoretically show that ${\rm DSA_2}$ enjoys {non-ergodic convergence properties}, i.e., the local minimizing sequence itself is convergent, a distinct feature that cannot be found in existing results. Specifically, we establish a convergence rate of $O(\frac{1}{\sqrt{t}})$ in terms of objective function error. Then, extensions are made to tackle distributed optimization problems with {coupled functional constraints} by combining ${\rm DSA_2}$ and dual decomposition. This is made possible by Lagrangian relaxation that transforms the coupling in constraints of the primal problem into that in cost functions of the dual, thus allowing us to solve the dual problem via ${\rm DSA_2}$. Both 
the dual objective error and the quadratic penalty for the coupled constraint are proved to converge at a rate of $O(\frac{1}{\sqrt{t}})$, and the primal objective error asymptotically vanishes.
Numerical experiments and comparisons are conducted to illustrate the advantage of the proposed algorithms and validate our theoretical findings. 
\end{abstract}

\end{frontmatter}

\section{Introduction}
This work considers large scale convex optimization problems that are defined over networks, and develops and analyzes distributed algorithms that are compatible with the communication constraints to solve them. Such optimization problems naturally arise in many engineering scenarios.
For example, problems such as estimation in sensor networks, distributed control of multi-agent systems, and resource allocation, can be formulated as distributed convex programs \cite{boyddistributed,nedicnetworktopology}. Advantages of distributed optimization over its centralized counterpart lie in that it offers a flexible and robust solution framework where only locally light computations and peer-to-peer communication are required to minimize a global objective function.

Due to their wide applications, distributed multi-agent decision making has been recently widely studied by researchers. In the literature, two types of distributed optimization problems are of particular interest, that is, optimization problems with \emph{coupled cost functions} or \emph{coupled constraints}.
They are essentially different in terms of the coupling sources that prevent decomposition of the original problem, thus making the design challenging. 

For \emph{optimization problems with coupled costs}, early distributed optimization algorithms can be found in the seminal work \cite{tsitsiklisdistributed}, where multiple processors cooperatively minimize a common objective function by conducting gradient-based local iterations and exchanging information with other agents asynchronously. 
Notable recent distributed optimization algorithms are reported in \cite{nedicachieving,yuanontheconvergence,xuconvergence,xuabregman,quharnessing,shiextra,shionthelinear,duchidualaveraging,hosseinionline}.
Technically speaking, their designs both contain the following two crucial steps. In the first step, one assigns local copies about the global decision variable to each node such that each node has a local version of the optimization variable to work with and imposes a consensus constraint on local estimates to guarantee the equivalence to the original problem. Then, a local iteration rule associated with an appropriate synchronization mechanism for updating local estimates of the global minimizer is designed.
 
Existing methods essentially differ from each other in terms of the design in the second step. 
Regarding the local update of global variables, the algorithms reported in \cite{nedicachieving,yuanontheconvergence,xuconvergence,xuabregman,quharnessing,shiextra,shionthelinear} update local estimates about the minimizer by using primal methods, which directly generate points in the feasible set that is contained in the primal space of variables.
Typical primal methods include the projected subgradient method, where the minimizing sequences are generated by shifting the test point along the opposite directions of subgradients and conducting Euclidean projections in an iterative way; 
for detailed introduction about primal and dual methods the readers are referred to the monograph \cite{nemirovskyproblem}.
In this class of methods, consensus among local estimates is usually enforced by averaging each agent's local estimate and the information received from its immediate neighbors at each round based on certain weight matrices, e.g., doubly stochastic matrices.
For unconstrained optimization problems with smooth objective functions, recent works make use of the dynamic average consensus scheme \cite{zhudiscretetime} to track the gradient of the overall objective, and consecutively move the local estimates along the opposite directions of the approximated overall gradients to achieve minimization.
Although the schemes developed in \cite{quharnessing,nedicachieving,xuconvergence} share similarities, the methodologies for establishing convergence properties are different from each other. For example, the authors in \cite{nedicachieving} and \cite{quharnessing} develop convergence results based on the small-gain theorem and linear system inequalities, respectively. 
These schemes provably enjoy faster convergence rate and exact consensual minimization; please see \cite{quharnessing} for more details.
There are also some distributed optimization algorithms available in the literature \cite{duchidualaveraging,hosseinionline,shahrampourdistributed} where the local iteration rule works in the dual space, e.g., dual averaging \cite{nesterovprimal}. It is shown in \cite{duchidualaveraging} that minimizing the dual model of the objective function can alleviate some technical difficulties caused by the projection step used in primal methods. 


It is worth noting that all the aforementioned algorithms may not be able to generate \emph{a convergent sequence of test points}. Indeed, they only guarantee convergence of the objective function values at the {running average} along the local minimizing sequence, i.e., ergodic convergence properties. This essentially allows undesired jumps of the objective function values at some iterations, possibly threatening the stability of the distributed system. In centralized optimization, this problem may be mitigated by further considering the best test point achieved so far. This procedure, however, may be not implementable in distributed scenarios since the quality test of certain points requires knowledge of the global objective. 

In another line of research, \emph{optimization problems with coupled constraints} have also recently been extensively studied in the literature \cite{nesterovdualsubgradient,nedicapproximate,falsonedualdecomposition,simonettoprimalrecovery,mateosdistributed,wangdistributedMPC,chatzipanagiotisanaugmented,notarnicolaconstraint,shersonontheduality,liangdistributed}. In this class of problems, each agent holds its own decision variable, objective function and constraints, and is coupled via global inequality constraints. A powerful methodology to this kind of problem is known as the Lagrangian relaxation that transforms the primal problem with shared constraints to the corresponding dual problem with coupled costs and solves it; see \cite{nedicapproximate,nesterovdualsubgradient,conejodecomposition} for more details. 
It is worth mentioning that, in such a framework having the optimal dual variable does not necessarily give us an optimal primal variable as the dual objective is generally nonsmooth at the optimal dual point, and thus nontrivial primal recovery schemes are required \cite{gustavssonprimal}.
Examples of this method include dual decomposition and augmented Lagrangian methods (also known as the method of multipliers) \cite{chatzipanagiotisanaugmented}.
Note that standard dual decomposition \cite{nedicapproximate,nesterovdualsubgradient} requires a fusion center that is able to communicate with all other agents to collect necessary gradient information of the dual objective function.
To enable fully distributed implementation, the work in \cite{wangdistributedMPC} forms a double-loop algorithm that combines the accelerated gradient method and a finite time consensus scheme to tackle the dual problem. 
{ Notably, the authors in \cite{liangdistributed} theoretically validate the use of a constant stepsize for the case with the objectives and the functions that characterize the coupled constraint being smooth.}
In a nonsmooth scenario, recent work in \cite{notarnicolaconstraint} properly relaxes the constraint-coupled problem and explores the duality principle twice to design a distributed iteration scheme. It is shown in \cite{notarnicolaconstraint} that the primal variable converges without any averaging steps.
Alternatively, the authors in \cite{simonettoprimalrecovery,falsonedualdecomposition,mateosdistributed} resort to the consensus-based distributed subgradient methods to solve the dual problem. To be specific, the authors in \cite{shersonontheduality} propose to use the alternating method of multipliers (ADMM) and the primal-dual method of multipliers (PDMM) to solve the dual problem; however they do not present convergence results. 
The work in \cite{simonettoprimalrecovery} establishes convergence rate for constant stepsizes and \cite{mateosdistributed} for decaying stepsizes both in virtue of the assumption that a Slater point exists and is known to all agents, which is somewhat restrictive although a construction process of a Slater point requiring extra negotiation between agents is provided in \cite{mateosdistributed}.
The framework considered in \cite{falsonedualdecomposition} relaxes this assumption, but the requirement on the stepsize is more restrictive, i.e., square summable stepsizes, and the convergence rate is missing. 
It is worth mentioning that the consensus-based distributed subgradient methods used in \cite{simonettoprimalrecovery,mateosdistributed,falsonedualdecomposition} for solving the dual problem cannot generate a convergent sequence of dual variables but their running averages or the best achieved dual variables, as we explained above.
As a consequence, the obtained dual variable sequence that plays an important role in allocating the coupled constraint resources does not necessarily stabilize the multi-agent system.

Note that the constraint-coupled optimization problems can be converted to problems with coupled costs by augmenting the local variable such that each agent becomes interested in a copy of the global variable. In doing so, the algorithms in \cite{nedicconstrainedconsensus,duchidualaveraging,shionthelinear} can be applied but with an increased communication and computation load.


This paper mainly contributes to this area in two aspects.
\begin{itemize}
\item The first contribution of this work is to provide a distributed subgradient method with double averaging (abbreviated as ${\rm DSA_2}$) for nonsmooth optimization that enjoys \emph{non-ergodic convergence properties}, i.e, the local minimizing sequence itself is convergent. The developed method is based on the centralized {subgradient method with double averaging} (${\rm SA_2}$) recently developed in \cite{nesterovquasimonotone,nesterovdualsubgradient}. However, the methodology for establishing convergence properties is significantly different from that in \cite{nesterovquasimonotone}, since in consensus-based distributed optimization one should carefully handle the inexact subgradient information and quantify the network effect caused by distributed implementation. Compared to existing distributed dual methods, e.g., distributed mirror descent \cite{shahrampourdistributed} and distributed dual averaging \cite{duchidualaveraging}, we further introduce an averaging step to the distributed minimization scheme and theoretically show that it is this extra averaging step that makes the sequence of local test points convergent. Since dual methods require a linear model to minimize at each round, the dynamic average consensus scheme is recruited to track the overall gradient as in \cite{quharnessing,xuconvergence} such that each agent maintains a local estimate of the global gradient to form the linear approximation of the global objective.
We establish an $O(\frac{1}{\sqrt{t}})$ convergence rate for the proposed strategy, which is known as the best achievable rate of convergence for subgradient methods. See Table \ref{overview} for a detailed comparison between the proposed algorithm and existing results.

\item
Extensions are made to solve large scale optimization problems with coupled functional constraints by combining ${\rm DSA_2}$ and dual decomposition.
By Lagrangian relaxation, the coupling in constraints in the primal problem is first transformed into that in objective functions of the dual problem.
Then a primal-dual sequence is constructed by solving the dual problem via ${\rm DSA_2}$ and using local estimates of the optimal dual to derive the corresponding primal variable.
A feature of this strategy is that agents only negotiate on dual variables but do not exchange information about local objective functions, constraints, and their optimal decisions, which can effectively help secure privacy among agents. 
We theoretically show that both the dual objective error and the quadratic penalty for the coupled constraint admit $O(\frac{1}{\sqrt{t}})$ upper bounds, and the primal objective error vanishes asymptotically.
Numerical simulations and comparisons with state-of-the-art algorithms verify our theoretical findings.
 
\end{itemize}

\emph{Notation}: We denote by $\mathbb{R}$ the set of real numbers and $\mathbb{R}^m$ the $m$-dimensional Euclidean space. 
In this space, we let $\lVert \cdot \rVert_p$ denote the $l_p$-norm operator, $\lVert \cdot \rVert_*$ the dual norm of $\lVert \cdot \rVert$, and $\langle \cdot,\cdot \rangle $ the inner product of two vectors. $0_m\in\mathbb{R}^m$ represents the vector of all zeros, and $\mathbf{1}$ stands for an $m$-dimensional all one column vector.
Notation `$\geq$' is element-wise when applied to vectors.
For a column vector $x\in\mathbb{R}^m$, 
$x(i)$ denotes the $i$th element of vector $x$.
$\triangle_m =\{ x\in\mathbb{R}^m \lvert x\geq 0_m,\sum_{i=1}^{m}x(i)=1 \}$ represents the $m$-dimensional probability simplex.
Given an $m\times m$ matrix $A$, we denote its singular values by $\sigma_1(A)\geq \sigma_2(A)\geq \cdots \geq \sigma_m(A) \geq 0$. 
A sequence $\{x_t\}_{t\geq 0}$ is said to have non-ergodic (ergodic) convergence rate $O(\cdot)$ if the rate is evaluated at the test point itself $x_t$ (a supporting running sequence $y_t=\frac{1}{t+1}\sum_{k=0}^{t}x_k$).

%
%
%

\begin{table*}
		\caption{An overview of existing distributed optimization algorithms. }
		\label{overview}
	\begin{tabular}{|c|c|c|c|c|c|c|c|c|c|c|}
		\hline
		\hline

		\multirow{2}{*}{Algorithms} &
		
	    \multirow{2}{*}{${\rm DSA_2}$}&
		
		\multicolumn{4}{|c|}{\bf{Coupled costs}} &
		
		\multicolumn{5}{|c|}{\bf{Coupled constraints}} \\

		\cline{3-11}   
		
	&	& \cite{duchidualaveraging} & \cite{nedicconstrainedconsensus} & \cite{quharnessing} & \cite{shiextra}& \cite{nedicapproximate} & \cite{nesterovdualsubgradient} & \cite{mateosdistributed}& \cite{falsonedualdecomposition}& {\cite{notarnicolaconstraint}}  \\
		
		\hline

		\multirow{2}{*}{Assumptions} & \multicolumn{3}{|c|}{\multirow{2}{*}{Convex, Constrained,}} & \multicolumn{2}{|c|}{\multirow{2}{*}{Convex, Smooth,}}  &  \multicolumn{2}{|c|}{\multirow{2}{*}{Convex, Constrained}} & \multicolumn{3}{|c|}{\multirow{2}{*}{Convex,}} \\
		
		&\multicolumn{3}{|c|}{Bounded subgradient}  & \multicolumn{2}{|c|}{Unconstrained}   & \multicolumn{2}{|c|}{Fusion center}  &  \multicolumn{3}{|c|}{Constrained}\\
		
		\hline
		Exactness & 	\multicolumn{5}{|c|}{Yes}& No & \multicolumn{4}{|c|}{Yes} \\
		
		\hline
		Iteration rule & 	\multicolumn{2}{|c|}{Dual methods} &	\multicolumn{4}{|c|}{Primal methods}  & Dual methods & \multicolumn{3}{|c|}{Primal methods} \\
		
		
		
		\hline
		\multirow{3}{*}{Convergence}&  \multicolumn{4}{|c|}{Objective error} & {Fixed point residual} & \multicolumn{5}{|c|}{Objective error} \\
		\cline{2-11}
		&\multicolumn{1}{|c|}{\bf{Non-ergodic}}&\multicolumn{4}{|c|}{Ergodic}    & \multirow{2}{*}{$O(\frac{1}{t})$}  &\multicolumn{2}{|c|}{\multirow{2}{*}{$O(\frac{1}{\sqrt{t}})$}} & \multicolumn{2}{|c|}{\multirow{2}{*}{N/A}} \\
		\cline{2-6} 
		&\multicolumn{1}{|c|}{$O(\frac{1}{\sqrt{t}})$}&\multicolumn{2}{|c|}{$O(\frac{\log(t)}{\sqrt{t}})$}  & \multicolumn{2}{|c|}{$O(\frac{1}{t})$} & & \multicolumn{2}{|c|}{}  & \multicolumn{2}{|c|}{} \\
		\hline
		
		\multirow{2}{*}{Remarks} & 	\multicolumn{10}{|c|}{{\multirow{2}{*}{Some distributed optimization algorithms based on augmented Lagrangian }}} \\
			&\multicolumn{10}{|c|}{method{ \cite{shionthelinear,shersonontheduality}} and operator splitting \cite{xuabregman} are not included due to limited space.}\\
		\hline

	\end{tabular}
	
\end{table*}

\section{Problem Statement and Preliminaries}
This section formally presents the distributed optimization problem and preliminaries about the centralized ${\rm SA_2}$ \cite{nesterovquasimonotone}. 

\subsection{Problem statement}

We consider a problem where $n$ agents connected via a network manage to collaboratively solve the following constrained optimization problem:
\begin{equation} \label{original_problem}
\min_{x\in \mathcal{X}} f(x) =\frac{1}{n}\sum_{i=1}^{n}f_i(x)
\end{equation}
where ${x}\in\mathbb{R}^m$ denotes the decision variable and $\mathcal{X}\subseteq \mathbb{R}^m$ the common closed convex constraint set. 
Throughout this paper, we assume without loss of generality that $0_m\in\mathcal{X}$, since it can always be met by translating $\mathcal{X}$.
Each function $f_i: \mathcal{X}\rightarrow\mathbb{R}$ that is convex and possibly nonsmooth represents the local objective privately known to agent $i$. Suppose that problem in \eqref{original_problem} admits at least one optimal solution. We denote by ${x}^*$ one of the minimizers and $f(x^*)$ the minimal function value. For function $f_i$, we denote by $\triangledown f_i(x)$ its arbitrary subgradient at $x\in\mathcal{X}$ that satisfies
\begin{equation*}
f_i(y)\geq f_i(x)+\langle\triangledown f_i(x),y-x \rangle, \forall y\in\mathcal{X}.
\end{equation*}


The following assumption is made for the objective function.
\begin{assumption}\label{Lipschitz_continuity}
	Each function $f_i$ is $L$-Lipschitz with respect to some norm $\lVert \cdot \rVert$, i.e.,
	\begin{equation}
	\lvert f_i(x)- f_i(y)  \rvert \leq L \lVert x-y\rVert, \forall x,y\in\mathcal{X}.
	\end{equation}
\end{assumption}
It is worth to mention that the Lipschitz continuity assumed in Assumption \ref{Lipschitz_continuity} holds for many functions, e.g., any convex function on a closed domain or polyhedral function on an arbitrary domain.
A consequence of this assumption is that we have all the subgradients of $f_i(x),\forall x\in\mathcal{X}$ bounded in the dual norm{ \cite{duchidualaveraging}}, i.e., 
\begin{equation*}
\lVert \triangledown f_i(x)\rVert_* \leq L.
\end{equation*} 

The communication network that connects the multi-agent system is modeled by an undirected and simple graph $\mathcal{G}=(\mathcal{V},\mathcal{E})$, where $\mathcal{V}=\{1,\cdots,n\}$ denotes the set of agents and $\mathcal{E}\subseteq \mathcal{V}\times\mathcal {V}$ the set of edges that correspond to the communication channels between agents. 
Note that this general graph imposes communication constraints on the working agents, that is, each agent can only communicate with its neighboring agents $ j\in \mathcal{N}_i=\{j\in \mathcal{V}|(i,j)\in \mathcal{E} \}$. 


\subsection{Subgradient method with double averaging}

This subsection briefly reviews the ${\rm SA_2}$, based on which the proposed algorithms are developed. We begin by introducing a {prox-function} $d:\mathcal{X}\rightarrow\mathbb{R}$ that enjoys the following properties.
\begin{assumption}\label{prox_function}
	1)  $d(x)\geq 0, \forall x\in\mathcal{X}$ and $d(0_m)=0$;
	2)  $d(x)$ is $1$-strongly convex on $\mathcal{X}$ with respect to the same norm as in Assumption \ref{Lipschitz_continuity}, i.e.,
	\begin{equation}
	d(y)\geq d(x)+\langle\triangledown d(x), y-x \rangle+\frac{1}{2}\lVert y-x\rVert^2, \forall x,y\in\mathcal{X}.
	\end{equation}
\end{assumption} 
We remark that this assumption is standard in the sense that it can be easily achieved by a large group of functions. For instance, the quadratic function $d(x)=\frac{1}{2}\lVert x\rVert^2_2$ satisfies $d(0_m)=0$ and is $1$-strongly convex with respect to the $l_2$-norm, { and the entropic function $
d(x)=\sum_{i=1}^{m}x(i)\log x(i) -x(i)
$
is strongly convex with respect to the $l_1$-norm for $x$ in the $m$-dimensional probability simplex $\triangle_m$.}

${\rm SA_2}$ generates sequences of estimates about the minimizer and the corresponding subgradient, i.e., $\{x_t\}_{t\geq 0}$ and $\{\frac{1}{t+1}\sum_{k=0}^{t}\triangledown f(x_k)\}_{t\geq 0}$,
in an iterative way.
In particular, the algorithm at each time stamp $t$ performs the following iteration
\begin{subequations}\label{centralized_DA_2}
	\begin{align}
	\hat{x}_{t+1}&=\arg \min_{x\in\mathcal{X}}\Big\{\big\langle \sum_{k=0}^{t}\triangledown f(x_k), x \big\rangle +\gamma_td(x) \Big\}\label{local_optimization_centralized} \\
	x_{t+1} &= \frac{t+1}{t+2}{x}_t+\frac{1}{t+2}\hat{x}_{t+1} =\frac{1}{t+2}\sum_{k=0}^{t+1}\hat{x}_{k}, \label{extra_averaging_step}
	\end{align}
\end{subequations}
where $\gamma_t$ is a non-decreasing sequence of positive parameters. This scheme shares similarities with other dual methods \cite{duchidualaveraging,nesterovprimal} where the calculation of test points involves iteratively minimizing an averaged linear approximation of the objective function $f$. 
However, directly minimizing a linear model of the objective may lead to oscillation. Thus in \eqref{local_optimization_centralized} a sum of the linear model and a weighted proximal function $d$ is minimized.
The averaging step in \eqref{extra_averaging_step} is a feature that cannot be found in other dual methods. Due to this feature, ${\rm SA_2}$ is able to produce a convergent minimizing sequence \cite{nesterovquasimonotone}.

\section{Distributed Subgradient Algorithm with Double Averaging}
In this section, we develop ${\rm DSA_2}$ and show its connections with some existing results.
\subsection{Development of ${ DSA_2}$}
Recall that our objective is to minimize the composite function in \eqref{original_problem} with distributed computations being conducted at each vertex of a connected graph. 
A classic technique to fulfill this task is to reformulate this problem as a consensus problem. That is, one assigns local copies of the global decision variable $x$ and the subgradient $\triangledown f(x)$ evaluated at the corresponding $x$, i.e., $x_{i}$ and $s_{i}$, to each agent $i$, and encodes agreement constraints on local estimates, i.e., $x_i=x_j,\forall i,j\in\mathcal{V}$, to ensure the equivalence to the original optimization problem. 
In doing so, each agent has local versions of the global variable and the corresponding subgradient information to operate with.

We observe from \eqref{centralized_DA_2} that $\hat{x}_{t}$ directly depends on the subgradient accumulated over time and contributes to the local update of $x_{t}$. 
This essentially implies that, to mimic the centralized minimization, the mechanism of updating $s_{i,t}$ to generate an accurate local estimate of the global subgradient is crucial. 
Further, it may be intuitive to expect that if the disagreement between $s_{i,t}$ and the true subgradient converges then minimization of the global objective can be achieved. 

To exploit this feature, we in this work employ the dynamic average consensus scheme \cite{zhudiscretetime} to track the subgradient of the aggregate objective function by using the local subgradient and the information gathered from immediate neighbors. More specifically, each agent properly weights the collected information at each iteration to generate an estimate of the global subgradient. 
To model this process, we assign 
a positive weight $p_{ij}$ to each communication link $(i,j)\in\mathcal{E}$ and leave $p_{ij}=0$ for other $(i,j)$ pairs.
We make the following standard assumption for the graph and the weight matrix $P=[p_{ij}]$.
\begin{assumption}\label{assumption_weight_matrix}
	1) The graph $\mathcal{G}$ is connected; 2) $P$ has a strictly positive diagonal, i.e., $p_{ii}>0$;
	3) $P$ is doubly stochastic, i.e., $P\mathbf{1}=\mathbf{1}$ and $\mathbf{1}^{\mathrm{T}}P=\mathbf{1}^{\mathrm{T}}$.
\end{assumption}
Note that Assumption \ref{assumption_weight_matrix} guarantees $\sigma_2(P)<1$, which can be illustrated as follows. Assumptions \ref{assumption_weight_matrix}-1 and \ref{assumption_weight_matrix}-2 make the matrix $P^{\mathrm{T}}P$ irreducible and primitive, respectively. This fact together with Assumption \ref{assumption_weight_matrix}-3 gives that $P^{\mathrm{T}}P$ has a unique Perron-Frobenius eigenvalue which is $1$, meaning that $\sigma_2(P)<1$. 

We now are equipped to present the subgradient tracking scheme{ \cite{varagnolonewton}}:
\begin{equation}
s_{i,t+1}=\sum_{j\in\mathcal{N}_i\cup\{i\}}p_{ij}s_{j,t}+\triangledown f_i(x_{i,t+1})-\triangledown f_i(x_{i,t}) . \label{dynamic_average_consensus}
\end{equation}
Denote $
\overline{s}_t=\frac{1}{n}\sum_{i=1}^{n}{s}_{i,t}$, and
$g_t=\frac{1}{n}\sum_{i=1}^{n}\triangledown f_i(x_{i,t})
$.
A lemma known as the conservation property is recalled { (Lemma 3 in \cite{xuconvergence})}.
\newtheorem{lemma}{Lemma}
\begin{lemma} \label{conservation_property}
	 If $s_{i,0}=\triangledown f_i(x_{i,0}),i\in\mathcal{V}$, then
$
	\overline{s}_{t+1}=g_{t+1}.
$
\end{lemma} 

Then, each agent is able to perform the following 
\begin{subequations}\label{local_iteration}
	\begin{align}
	\hat{x}_{i,t+1}&=\arg \min_{x\in\mathcal{X}}\Big\{\big\langle \sum_{k=0}^{t}s_{i,k}, x \big\rangle +\gamma_td(x)\Big \} \label{local_optimization}\\
	x_{i,t+1} &= \frac{t+1}{t+2}{x}_{i,t}+\frac{1}{t+2}\hat{x}_{i,t+1}=\frac{1}{t+2}\sum_{k=0}^{t+1}\hat{x}_{i,k}, \label{local_averaging}
	\end{align}
\end{subequations}
where $x_{i,t}$ denotes the estimated variable maintained by agent $i$ at time stamp $t$. It is worth mentioning that the difference between Eqs. \eqref{local_optimization_centralized} and \eqref{local_optimization} is that \eqref{local_optimization_centralized} uses exactly the accumulated subgradient of the global objective, i.e., $\triangledown f(x_t)$, while \eqref{local_optimization} uses an estimated one, i.e., $s_{i,t}$, due to the incomplete knowledge of each agent about the global objective. 

The proposed ${\rm DSA_2}$ is detailed in Algorithm \ref{DSA_2}.
\begin{algorithm}
	\begin{algorithmic}[1]
		\caption{${\rm DSA_2}$}
		\label{DSA_2}
		\STATE  Set $t=0, s_{i,0}=\triangledown f_i(x_{i,0})$, choose a non-decreasing sequence of positive parameters $\{\gamma_t\}_{t\geq 0}$.		
		\WHILE  { { Convergence is not reached}}
		\FOR{Each agent $i\in\mathcal{V}$ (in parallel)}  
		\STATE  Receive $s_{j,t}, \forall j\in\mathcal{N}_{i}$;
		\STATE Perform local computation in \eqref{local_iteration} and \eqref{dynamic_average_consensus};
		\STATE Broadcast $s_{i,t+1}$ to $j\in\mathcal{N}_i$;
		\ENDFOR 
		\STATE Set $t = t+1$.
		\ENDWHILE
	\end{algorithmic}
\end{algorithm}

\newtheorem{remark}{Remark}
\begin{remark}
It is worth mentioning that in most of the existing results, the running average is merely used as a supporting sequence for convergence but not involved in the subproblems at each iteration. That is, they have a recursion rule similar to the following{ \cite{duchidualaveraging}}
\begin{equation*}
\begin{split}
\hat{x}_{i,t+1}&=\arg \min_{x\in\mathcal{X}}\Big\{\big\langle \sum_{k=0}^{t}s_{i,k}, x \big\rangle +\gamma_td(x) \Big\} \\
s_{i,t+1}&=\sum_{j=1}^{n}p_{ij}s_{j,t}+\triangledown f_i(\hat{x}_{i,t+1})-\triangledown f_i(\hat{x}_{i,t}),
\end{split}
\end{equation*}
and establish convergence for the supporting running sequence $	x_{i,t+1}= \frac{t+1}{t+2}{x}_{i,t}+\frac{1}{t+2}\hat{x}_{i,t+1}=\frac{1}{t+2}\sum_{k=0}^{t+1}\hat{x}_{i,k}$. 
However, sometimes it is more desirable to have a convergent sequence of test points, taking the decentralized dual Lagrangian problem for example where the local test point is further used to coordinate subproblems. In this work, we show that, upon using the subgradient information evaluated at the running sequence $x_{i,t+1}$, i.e., Eq. \eqref{dynamic_average_consensus}, a local convergent sequence of test points can be obtained.
\end{remark}

\subsection{Relation to existing results}
As mentioned in Introduction, existing distributed optimization methods contain two important parts, that is, an update rule for local estimates about the global variable and a properly designed consensus mechanism. In this subsection, we compare ${\rm DSA_2}$ with existing works in terms of these two.

Regarding the local update of global variables, the proposed algorithm works in the dual space, i.e., each agent at each round generates a linear model of the local objective function and maps it back into the primal space by performing a projection-like operation, i.e., \eqref{local_optimization}.
Notable works including \cite{duchidualaveraging,hosseinionline} enjoy this feature too. The difference between them and the proposed one is that we further introduce an averaging step \eqref{local_averaging}
that allows us to generate {a convergent sequence of test points}.
Thanks to this distinct feature, the proposed algorithm further lends itself to consensus-based dual decomposition that is known as a powerful methodology to deal with large scale \emph{constraint-coupled} optimization problems, as we will show later.  
We now turn to compare the consensus mechanisms. 
Define $z_{i,t}=\sum_{k=0}^{t}s_{i,k}$. It can be seen that if $s_{i,0}=\triangledown f_i(x_{i,0})$
then the subgradient tracking scheme \eqref{dynamic_average_consensus} reduces to
\begin{equation*}
z_{i,t+1}=\sum_{j\in\mathcal{N}_i\cup\{i\}}p_{ij}z_{j,t}+\triangledown f_i(x_{i,t+1}),
\end{equation*}
which is the subgradient update rule in \cite{duchidualaveraging}.

In another line of research, local estimates about the minimizer are updated by using primal methods which directly generate a sequence of points in the feasible set that is contained in the primal space of variables \cite{nedicconstrainedconsensus,simonettoprimalrecovery}. Specifically, each agent at each round conducts the following
\begin{equation}\label{primal_method}
x_{i,t+1}= \mathcal{P}_{\mathcal{X}}\Big[ \sum_{j\in\mathcal{N}_i\cup\{i\}}p_{ij}x_{j,t}-\frac{1}{\gamma_{t}}\triangledown f_i(x_{i,t}) \Big],
\end{equation}
where $\mathcal{P}_{\mathcal{X}}\left[ y \right] =\arg \min_{x\in\mathcal{X}}\|y-x \|_2$ denotes the Euclidean projection of a vector $y$ on the set $\mathcal{X}$.
Note that there is a consensus seeking step before the projection in \eqref{primal_method}. For unconstrained optimization problems with smooth objective functions, recent works \cite{xuconvergence,quharnessing,nedicachieving} make use of the gradient tracking scheme \eqref{dynamic_average_consensus} and modify the rule to update local variables in \eqref{primal_method} as
\begin{equation*}
x_{i,t+1}=  \sum_{j\in\mathcal{N}_i\cup\{i\}}p_{ij}x_{j,t}-\frac{1}{\gamma}s_{i,t}, 
\end{equation*}
where $\frac{1}{\gamma}$ denotes the constant stepsize. 

We remark that, to the best of our knowledge, the proposed method is the first distributed optimization method being able to establish non-ergodic convergence rate in terms of objective function error.


\section{Convergence Properties Analysis}
\subsection{Basic convergence analysis}
In this subsection, we state the basic convergence results that reveal how the local estimate $x_{i}$ approaches the minimizer with the help of the global subgradient tracking scheme \eqref{dynamic_average_consensus}, highlighting the network effect due to distributed implementation.


Motivated by the literature regarding consensus-based distributed optimization, we set up an auxiliary sequence $\{y_t\}_{t\geq 0}$ that makes use of the {averaged subgradient} $g_t$ and conducts the following recursion
\begin{equation*}
	\begin{split}
	\hat{y}_0&=y_0 = 0_m\\
	\hat{y}_{t+1}&=\arg \min_{x\in\mathcal{X}}\Big\{\big\langle \sum_{k=0}^{t}g_{k}, x \big\rangle +\gamma_td(x) \Big\} \\
	y_{t+1} &= \frac{t+1}{t+2}{y}_{t}+\frac{1}{t+2}\hat{y}_{t+1}.
	\end{split}
\end{equation*}
Note that we only impose specific initial conditions for the sequence $\{y_t\}_{t\geq 0}$; the initial guess $x_{i,0}$ for each agent can be arbitrary. 

We present a slightly modified result in dual averaging { (Theorem 2 in \cite{nesterovprimal}, Lemma 3 in \cite{duchidualaveraging})}. Note that we let by convention that $\gamma_{-1}=\gamma_0$.

\begin{lemma} \label{auxiliary_sequence}
	Suppose Assumptions \ref{Lipschitz_continuity}-\ref{prox_function} hold true. For any non-decreasing sequence $\{\gamma_t\}_{t\geq 0}$ of positive parameters, and $x\in\mathcal{X}$, we have
	\begin{equation*}
	\sum_{k=0}^{t}\langle g_k, \hat{y}_k-x \rangle\leq \frac{1}{2} \sum_{k=0}^{t}\frac{1}{\gamma_{k-1}}\lVert g_k\rVert_*^2+\gamma_t d(x).
	\end{equation*} 
\end{lemma}

To establish relations between the local estimate $\{x_t\}_{t\geq 0}$ and $\{y_t\}_{t\geq 0}$, we recall the following standard result in convex analysis { (Lemma 1 in \cite{nesterovprimal})}.
\begin{lemma} \label{gamma_continuity}
	For any $u,v\in\mathbb{R}^m$ and $\gamma>0$, we have
	\begin{equation*}
	\begin{split}
	&\big\lVert \arg \min_{x\in \mathcal{X}}\big\{ \langle u,x \rangle+\gamma d(x) \big\} -\arg \min_{x\in \mathcal{X}}\big\{ \langle v,x \rangle+\gamma d(x) \big\} \big\rVert \\
	& \leq \frac{1}{\gamma} \lVert u-v \rVert_*.
	\end{split}
	\end{equation*}
\end{lemma}

\begin{lemma}\label{basic_convergence_thm}
	Suppose Assumptions \ref{Lipschitz_continuity}-\ref{prox_function} hold true.
	Let the sequences $\{x_{i,t}\}_{t\geq 0}$ and $\{s_{i,t}\}_{t\geq 0}$ be generated by Algorithm \ref{DSA_2}. For any $x\in\mathcal{X}$, we have
\begin{equation}\label{basic_convergence}
\begin{split}
&f(x_{i,t})-f(x)
\leq  \frac{L}{t+1}\times\\
& \sum_{k=0}^{t} \frac{1}{\gamma_{k-1}} \Big(\big\lVert \sum_{l=0}^{k-1}( s_{i,l}-g_l ) \big\rVert_* +\frac{2}{n}\sum_{j=1}^{n}\big\lVert \sum_{l=0}^{k-1}( s_{j,l}-g_l ) \big\rVert_* \Big) \\
& + \frac{1}{2(t+1)
}  \sum_{k=0}^{t}\frac{1}{\gamma_{k-1}}\lVert g_k\rVert_*^2+\frac{1}{t+1} \gamma_t d(x).
\end{split}
\end{equation}
\end{lemma}
\begin{pf}
For any $x\in\mathcal{X}$, we consider
\begin{equation}\label{ergodic_and_non}
\begin{split}
(t+1)&\sum_{j=1}^{n}\big(f_j(x_{j,t})-f_j(x)\big) = (t+1) \sum_{j=1}^{n}f_j(x_{j,t})\\
&-\sum_{k=0}^{t}\sum_{j=1}^{n}f_j(x_{j,k}) +\sum_{k=0}^{t}\sum_{j=1}^{n}\big(f_j(x_{j,k})-f_j(x)\big).
\end{split}
\end{equation}
By convexity of $f_i$, we have
\begin{equation}\label{convexity_1}
\begin{split}
&(t+1) \sum_{j=1}^{n}f_j(x_{j,t})-\sum_{k=0}^{t}\sum_{j=1}^{n}f_j(x_{j,k}) \\
=&  t \sum_{j=1}^{n}f_j(x_{j,t})-\sum_{k=0}^{t-1}\sum_{j=1}^{n}f_j(x_{j,k}) \\
=& \sum_{j=1}^{n}\sum_{k=1}^{t} k \big(f_j(x_{j,k})-f_j(x_{j,k-1})\big) \\
\leq & \sum_{j=1}^{n}\sum_{k=1}^{t}k\langle\triangledown f_j(x_{j,k}), x_{j,k}-x_{j,k-1} \rangle 
\end{split}
\end{equation}
and
\begin{equation}\label{conveixty_2}
\begin{split}
\sum_{j=1}^{n}\big(f_j(x_{j,k})-f_j(x)\big)\leq   \sum_{j=1}^{n}\langle \triangledown f_j(x_{j,k}), x_{j,k}-x \rangle.
\end{split}
\end{equation}
Plugging inequalities \eqref{convexity_1} and \eqref{conveixty_2} into \eqref{ergodic_and_non} yields
\begin{equation*}
\begin{split}
&(t+1)\sum_{j=1}^{n}\big(f_j(x_{j,t})-f_j(x)\big) \\
 \leq & \sum_{j=1}^{n}\Big( \sum_{k=1}^{t}  \langle \triangledown f_j(x_{j,k}), (k+1)x_{j,k}-kx_{j,k-1}-x \rangle \\
& + \langle \triangledown f_j(x_{j,0}), x_{j,0}-x \rangle \Big),
\end{split}
\end{equation*}
which in conjunction with an equivalent expression of \eqref{local_averaging}
\begin{equation*}
(t+1)x_{j,t}=tx_{j,t-1}+\hat{x}_{j,t}
\end{equation*}
gives rise to
\begin{equation*}
(t+1)\sum_{j =1}^{n}\big(f_j(x_{j,t})-f_j(x)\big) \leq \sum_{j=1}^{n} \sum_{k=0}^{t} \langle \triangledown f_j(x_{j,k}), \hat{x}_{j,k}-x \rangle .
\end{equation*}
To approach the desired result using Lemma \ref{auxiliary_sequence}, we shall rewrite the above inequality as 
\begin{equation} \label{crucial_inequality}
\begin{split}
&(t+1)\sum_{j=1}^{n}\big(f_j(x_{j,t})-f_j(x)\big)\\
 \leq & \sum_{j=1}^{n} \sum_{k=0}^{t} \Big( \langle \triangledown f_j(x_{j,k}), \hat{x}_{j,k}-\hat{y}_{k} \rangle+
 \langle \triangledown f_j(x_{j,k}), \hat{y}_{k}-x \rangle \Big) \\
 =& \sum_{k=0}^{t}  \Big( \sum_{j=1}^{n} \langle \triangledown f_j(x_{j,k}), \hat{x}_{j,k}-\hat{y}_{k} \rangle +  n\langle g_k, \hat{y}_{k}-x \rangle \Big).
\end{split}
\end{equation}
With this relation in mind, we now turn to consider
\begin{equation}\label{identical_variable}
\begin{split}
&f(x_{i,t})-f(x) = f(x_{i,t})-f(y_t)+f(y_t)-f(x) \\
\leq& \frac{1}{n}\Big(  \sum_{j=1}^{n} \big(f_j(y_t)-f_j(x_{j,t})\big) +\sum_{j=1}^{n}\big(f_j(x_{j,t})-f_j(x) \big)\Big)\\
&+  L\lVert x_{i,t}-y_t \rVert \\
\leq& L \Big(\lVert x_{i,t}-y_t \rVert +\frac{1}{n}\sum_{j=1}^{n}\lVert x_{j,t}-y_t \rVert \Big)\\ &+ \frac{1}{n} \sum_{j=1}^{n}\big(f_j(x_{j,t})-f_j(x) \big),
\end{split}
\end{equation}
where we use the $L$-Lipschitz continuity of $f_i$ to derive the first and second inequality. In light of \eqref{crucial_inequality}, we have
\begin{equation*}
\begin{split}
&f(x_{i,t})-f(x)\\
 &\leq  L \Big(\lVert x_{i,t}-y_t \rVert +\frac{1}{n}\sum_{j=1}^{n}\lVert x_{j,t}-y_t \rVert \Big) \\
& + \frac{1}{t+1} \sum_{k=0}^{t}\Big(\frac{1}{n} \sum_{j=1}^{n}\langle \triangledown f_j(x_{j,k}), \hat{x}_{j,k}-\hat{y}_k \rangle+\langle g_k, \hat{y}_k-x \rangle   \Big).
\end{split}
\end{equation*}
By the fact that
$
y_t= \frac{1}{t+1}\Big(y_0+\sum_{k=1}^{t}\hat{y}_k\Big)= \frac{1}{t+1}\sum_{k=0}^{t}\hat{y}_k
$,
we obtain
\begin{equation*}
\begin{split}
&f(x_{i,t})-f(x)\\
\leq & \frac{L}{t+1} \sum_{k=0}^{t}\Big(\lVert \hat{x}_{i,k}-\hat{y}_k \rVert +\frac{1}{n}\sum_{j=1}^{n}\lVert \hat{x}_{j,k}-\hat{y}_k \rVert \Big) \\
 +& \frac{1}{t+1} \sum_{k=0}^{t}\Big(\frac{1}{n} \sum_{j=1}^{n}\langle \triangledown f_j(x_{j,k}), \hat{x}_{j,k}-\hat{y}_k \rangle+\langle g_k, \hat{y}_k-x \rangle   \Big)\\
\leq  & \frac{L}{t+1} \sum_{k=0}^{t}\Big(\lVert \hat{x}_{i,k}-\hat{y}_k \rVert +\frac{2}{n}\sum_{j=1}^{n}\lVert \hat{x}_{j,k}-\hat{y}_k \rVert \Big ) \\
 & + \frac{1}{t+1}\sum_{k=0}^{t} \langle g_k, \hat{y}_k-x \rangle.
\end{split}
\end{equation*}
It follows from Lemma \ref{auxiliary_sequence} that
\begin{equation*}
\begin{split}
&f(x_{i,t})-f(x)\\
\leq & \frac{L}{t+1} \sum_{k=0}^{t}\Big(\lVert \hat{x}_{i,k}-\hat{y}_k \rVert +\frac{2}{n}\sum_{j=1}^{n}\lVert \hat{x}_{j,k}-\hat{y}_k \rVert \Big) \\
& + \frac{1}{2(t+1)
}  \sum_{k=0}^{t}\frac{1}{\gamma_{k-1}}\lVert g_k\rVert_*^2+\frac{1}{t+1} \gamma_t d(x).
\end{split}
\end{equation*}

Appealing to the $\frac{1}{\gamma}$-Lipschitz continuity of the projection operator in \eqref{local_optimization} (Lemma \ref{gamma_continuity}) allows us to obtain the desired result in \eqref{basic_convergence}. \QEDA 
\end{pf}

\begin{remark}
	Lemma \ref{basic_convergence_thm} highlights that, after $t$ steps of execution, the objective error $f(x_{i,t})-f(x^*)$ is bounded from above by a summation of four terms. The first two terms are due to the different estimates of the averaged subgradient of the global objective. The third and the fourth terms can be seen as the optimization error terms observed also in centralized nonsmooth optimization. The result suggests that, if the deviation $\big\lVert \sum_{l=0}^{k-1}( s_{i,l}-g_l ) \big\rVert_*$ is finite, i.e., $\lVert  s_{i,l}-g_l  \rVert_*$ decays fast enough, and $\gamma_{k}$ is properly chosen, then the objective error asymptotically converges to $0$. The assumptions for Lemma \ref{basic_convergence_thm} (Assumptions \ref{Lipschitz_continuity}-\ref{prox_function}) are not restrictive in the sense that the Lipschitz continuity holds for general convex functions defined on a closed domain and the requirements for the prox-function in Assumption \ref{prox_function} can also be easily met.
\end{remark}

{ 
\begin{remark}
	The proposed algorithm needs a synchronized decreasing $\frac{1}{\gamma_t}$ to ensure convergence. The choice is made for technical reasons. Specifically, it has to be made identical for each agent to validate the use of the $\frac{1}{\gamma}$-Lipschitz continuity of the projection operator in \eqref{local_optimization}, and to be decaying to make the objective error convergent in light of Eq. \eqref{basic_convergence}. This requirement can be satisfied by a synchronization step before execution of the algorithm. We note that the work in \cite{xuconvergence} relaxes this requirement, where an unconstrained distributed optimization problem is considered.
\end{remark}}

\subsection{Disagreement analysis}

\begin{lemma}\label{disagreement}
		Suppose Assumptions \ref{Lipschitz_continuity} and \ref{assumption_weight_matrix} hold. For the sequences $\{s_{i,t}\}_{t\geq 0}$ generated by the subgradient tracking scheme \eqref{dynamic_average_consensus}, we have
	\begin{equation}\label{disagreement_bound}
	\big\lVert \sum_{l=0}^{k}s_{i,l}-\sum_{l=0}^{k}g_{l}\big\rVert_*
	\leq \frac{\sqrt{n}L}{1-\sigma_2(P)}+2L.
	\end{equation}
\end{lemma}
\begin{pf}
	Define 
	\begin{equation*}
	\begin{split}
	{s}_t=\begin{bmatrix}
	s_{1,t}\\ \vdots \\ s_{n,t}
	\end{bmatrix},
	{\triangledown}_t=\begin{bmatrix}
	{\triangledown}f_1(x_{1,t})\\ \vdots \\{\triangledown}f_n(x_{n,t})
	\end{bmatrix}
	\end{split} 
	\end{equation*} and rewrite the dynamics in \eqref{dynamic_average_consensus} for all $i$ in a compact form as
	\begin{equation}  \label{gradient_error}
	s_{t+1}=Ps_t+\triangledown_{t+1}-\triangledown_{t}.
	\end{equation}
	Summing \eqref{gradient_error} over $t$ from $t=0$ to $k-1$ yields
	$
	\sum_{l=1}^{k}s_{l} =P\sum_{l=0}^{k-1}s_{l}-\triangledown_{0}+\triangledown_{k}
	$.
	Since $s_0=\triangledown_0$, we have
	\begin{equation*}
	\begin{split}
	&\sum_{l=0}^{k}s_l =P\sum_{l=0}^{k-1}s_{l}+\triangledown_{k}.
	\end{split}
	\end{equation*}
	To establish relation between $\sum_{l=0}^{k}s_l$ and the accumulated averaged subgradient, we subtract $\sum_{l=0}^{k}g_l$ on both sides and get
	\begin{equation*}
	\begin{split}
	\sum_{t=0}^{k}s_{t}-\mathbf{1}\sum_{t=0}^{k}g_{t} =&P\sum_{t=0}^{k-1}s_{t}-\mathbf{1}\sum_{t=0}^{k-1}g_{t}-\mathbf{1}g_{k}+\triangledown_{k}.
	\end{split}
	\end{equation*}
	By using  $\big(P-\frac{1}{n}\mathbf{1}\mathbf{1}^{\mathrm{T}}\big)(s_t-\mathbf{1}g_t)=Ps_t-\mathbf{1}g_t$, we obtain
	\begin{equation*}
	\begin{split}
	\sum_{l=0}^{k}s_{l}&-\mathbf{1}\sum_{l=0}^{k}g_{l}\\
	=&\big(P-\frac{1}{n}\mathbf{1}\mathbf{1}^{\mathrm{T}}\big)\Big(\sum_{l=0}^{k-1}s_{l}-\mathbf{1}\sum_{l=0}^{k-1}g_{l}\Big)-\mathbf{1}g_{k}+\triangledown_{k}.
	\end{split}
	\end{equation*}
	By recursion, we get
	\begin{equation*}
	\begin{split}
	&\sum_{l=0}^{k}s_{l}-\mathbf{1}\sum_{l=0}^{k}g_{l}\\ =&\sum_{l=0}^{k-1}\big(P-\frac{1}{n}\mathbf{1}\mathbf{1}^{\mathrm{T}}\big)^{k-l}(\triangledown_{l}-\mathbf{1}g_l)-\mathbf{1}g_{k}+\triangledown_{k} \\
	=&  \sum_{l=0}^{k-1}\big(P^{k-l}-\frac{1}{n}\mathbf{1}\mathbf{1}^{\mathrm{T}}\big)\triangledown_{l}-\mathbf{1}g_{k}+\triangledown_{k},
	\end{split}
	\end{equation*}
	implying that
	\begin{equation*}
	\begin{split}
	&\sum_{l=0}^{k}s_{i,l}-\sum_{l=0}^{k}g_{l}\\
	& = \sum_{l=0}^{k-1}\sum_{j=1}^{n}\big(P^{k-l}-\frac{1}{n}\mathbf{1}\mathbf{1}^{\mathrm{T}}\big)_{ij}\triangledown f_j(x_{j,l})-g_{k}+\triangledown f_i(x_{i,k}).
	\end{split}
	\end{equation*}
	Taking the dual norm on both sides gives rise to
	\begin{equation*}
	\begin{split}
	\big\lVert &\sum_{l=0}^{k}s_{i,l}-\sum_{l=0}^{k}g_{l}\big\rVert_*\\
	 =& \sum_{l=0}^{k-1}\sum_{j=1}^{n}\big\lvert\big(P^{k-l}-\frac{1}{n}\mathbf{1}\mathbf{1}^{\mathrm{T}}\big)_{ij}\big\rvert\lVert\triangledown f_j(x_{j,l})\rVert_*\\
	&+\lVert\triangledown f_i(x_{i,k})-g_{k}\rVert_* 
	\leq \sum_{l=0}^{k-1}\big\lVert P^{k-l}e_i-\frac{\mathbf{1}}{n} \big\rVert_1L+2L,
	\end{split}
	\end{equation*}
	where $e_i\in\mathbb{R}^m $ denotes the $i$-th standard basis vector.
	Recall that for a stochastic matrix $P$ \cite{duchidualaveraging,hosseinionline} one has
	$
	\big\lVert  P^{k-l}e-\frac{\mathbf{1}}{n} \big\rVert_1\leq \sqrt{n}\big\lVert  P^{k-l}e-\frac{\mathbf{1}}{n} \big\rVert_2 \leq \sigma_2(P)^{k-l}\sqrt{n}
	$
	for all $e\in\triangle_m$.
	The inequality in \eqref{disagreement_bound} then follows, thereby concluding the proof. \QEDA 
\end{pf}

\subsection{Non-ergodic convergence rate}
We now are in a position to establish the non-ergodic convergence rate.
\begin{thm} \label{convergence_rate}
	Suppose that $d(x^*)\leq R^2$ and $\gamma_{t}=\gamma \sqrt{t+1}$ where $\gamma> 0$, and Assumptions \ref{Lipschitz_continuity}-\ref{assumption_weight_matrix} hold. For the sequences $\{x_{i,t}\}_{t\geq 0}$ and $\{s_{i,t}\}_{t\geq 0}$ being generated by Algorithm \ref{DSA_2}, we have
	\begin{equation}\label{convergence_rate_bound}
	\begin{split}
	&f(x_{i,t})-f(x^*)\\
	\leq & \frac{1}{\sqrt{t+1}}\Big( \big( \frac{6L^2\sqrt{n}}{1-\sigma_2(P)}+13L^2\big)\frac{1}{\gamma}  +{\gamma R^2}\Big) .
	\end{split}
	\end{equation}
\end{thm}
\begin{pf}
By invoking Lemma \ref{disagreement} and boundedness of $\lVert g_k \rVert_*$, we can obtain from the result in Lemma \ref{basic_convergence_thm} that
\begin{equation*}
\begin{split}
&f(x_{i,t})-f(x)\\
\leq & \frac{3L^2}{t+1} \big(\frac{\sqrt{n}}{1-\sigma_2(P)}+2\big)  \sum_{k=0}^{t} \frac{1}{\gamma_{k-1}}  \\
& + \frac{1}{2(t+1)
}  \sum_{k=0}^{t}\frac{1}{\gamma_{k-1}}\lVert g_k\rVert_*^2+\frac{1}{t+1} \gamma_t d(x) \\
\leq &  \frac{3L^2}{t+1} \big(\frac{\sqrt{n}}{1-\sigma_2(P)}+2\big)  \sum_{k=0}^{t} \frac{1}{\gamma_{k-1}}  \\
& + \frac{L^2}{2(t+1)
}  \sum_{k=0}^{t}\frac{1}{\gamma_{k-1}}+\frac{1}{t+1} \gamma_t d(x).
 \end{split}
\end{equation*}
Due to the fact that
\begin{equation*}
\begin{split}
&\sum_{k=0}^{t}\frac{1}{\gamma_{k-1}} = \frac{1}{\gamma_0}+\sum_{k=0}^{t-1}\frac{1}{\gamma_{k}} = \frac{1}{\gamma}+\frac{1}{\gamma}\sum_{k=0}^{t-1}\frac{1}{\sqrt{k+1}} \\
& \leq \frac{2}{\gamma}\sqrt{t+1},
\end{split}
\end{equation*}
we get
\begin{equation*}
\begin{split}
&f(x_{i,t})-f(x)\\
\leq & \big( \frac{6L^2\sqrt{n}}{1-\sigma_2(P)}+13L^2\big)\frac{1}{\gamma\sqrt{t+1}}  +\frac{\gamma}{\sqrt{t+1}}d(x).
\end{split}
\end{equation*}

We arrive at the desired result in \eqref{convergence_rate_bound} by using the assumption that $d(x^*)\leq R^2$. \QEDA 
	\end{pf}


\begin{remark}
	In Theorem \ref{convergence_rate}, we establish a non-ergodic convergence rate $O(1/\sqrt{t})$ for nonsmooth objective functions in terms of the objective error. 
	The result is a little bit stronger than the ergodic one in most of the existing results. To see this, we note
	\begin{equation*}
	f(\frac{1}{t}\sum_{k=1}^{t}x_{i,k})-f(x^*)\leq \frac{1}{t}\sum_{k=1}^{t}\big(f(x_{i,k})-f(x^*)\big) \leq O(1/\sqrt{t})
	\end{equation*}
	by convexity of $f$. More importantly, this essentially makes the proposed method applicable to decentralized dual Lagrangian problems where the local test point is further used to coordinate subproblems, as we shall see in the next section.
\end{remark}

Note that the result in Theorem \ref{convergence_rate} allows us to further derive the optimal choice of $\gamma$, i.e.,
$
\gamma = \frac{L}{R}\sqrt{\frac{6\sqrt{n}}{1-\sigma_2(P)}+13},
$
with which the convergence rate result becomes
\begin{equation*}
\begin{split}
f(x_{i,t})-f(x^*)
\leq \frac{2RL}{\sqrt{t+1}} \sqrt{\frac{6\sqrt{n}}{1-\sigma_2(P)}+13}.
\end{split}
\end{equation*}
It is shown in \eqref{convergence_rate_bound} that the asymptotical bound depends on $\sigma_2(P)$. This fact can be explored to achieve a tighter bound by following the methods in \cite{xiaodistributed} to minimize the spectral norm of $P-\frac{1}{n}\mathbf{1}\mathbf{1}^{\mathrm{T}}$. For special graphs such as paths and cycles with a specific structure of weight matrix $P$, $\sigma_2(P)$ can be explicitly identified to make the bound tighter \cite{duchidualaveraging}.

\section{Extension to Constraint-Coupled Distributed Optimization}\label{section5}
In this section, we further develop a ${\rm DSA_2}$-based dual decomposition strategy to solve optimization problems with coupled functional constraints.

Consider the following minimization problem
\begin{equation} \label{coupled_constraints}
\begin{split}
\min_{\{x_i\in \mathcal{X}_i\}_{i=1}^n}  &\sum_{i=1}^{n}f_i(x_i) \\
 {\rm subject \,\ to} \quad &\sum_{i=1}^{n}h_i(x_i)\leq 0_m,
\end{split}
\end{equation}
where $\mathcal{X}_i$ are compact and convex sets, and $f_i: \mathcal{X}_i\rightarrow\mathbb{R}$ and $h_i: \mathcal{X}_i\rightarrow\mathbb{R}^m$ are closed and convex functions. 
We denote by $\{x_i^*\}_{i=1}^n$ one of the optimal solutions and $\sum_{i=1}^{n}f_i^*$ the minimal function value.

We can see from the above problem that the objective function and parts of the constraints, i.e., $\{x_i\in\mathcal{X}_i\}_{i=1}^n$, enjoy a separable structure, but the global constraint $\sum_{i=1}^{n}h_i(x_i)\leq 0_m$ cannot be trivially decomposed. 
One powerful methodology to solve this problem is to alternatively consider the corresponding dual Lagrangian problem.
In doing so, the coupling in constraints can be transformed into that in objective functions, thus allowing us to solve it via the proposed ${\rm DSA_2}$.

The Lagrangian of \eqref{coupled_constraints} is
\begin{equation*}
\begin{split}
&L(\{x_i\}_{i=1}^n,\lambda)\\
=&\sum_{i=1}^{n}L_i({x}_i,\lambda)=\sum_{i=1}^{n}\big(f_i(x_i)+\langle \lambda, h_i(x_i) \rangle \big),
\end{split}
\end{equation*} 
where $x_i\in\mathcal{X}_i$ and $\lambda\geq0_m$ represents the dual variable associated with the coupled constraint, and the dual Lagrangian problem is
\begin{equation*}
\max_{\lambda\geq0_m} \min_{\{x_i\in\mathcal{X}_i\}_{i=1}^n} \sum_{i=1}^{n}L_i(x_i,\lambda),
\end{equation*}
which is equivalent to 
\begin{equation*} 
 \min_{\lambda\geq0_m}\sum_{i=1}^{n}\psi_i(\lambda) =\min_{\lambda\geq0_m} \max_{\{x_i\in\mathcal{X}_i\}_{i=1}^n} -\sum_{i=1}^{n}L_i(x_i,\lambda).
\end{equation*}

It is worth mentioning that the dual Lagrangian problem has the same structure as that in \eqref{original_problem} treated in previous sections.
To see this, we define
$
\psi_i(\lambda)= \max_{x_i\in\mathcal{X}_i}-L_i(x_i,\lambda)
$
and rewrite the dual Lagrangian problem as
\begin{equation} \label{dual_Lagrangian}
\min_{\lambda\geq0_m} \sum_{i=1}^{n}\psi_i(\lambda).
\end{equation}
It is assumed that the dual Lagrangian problem is solvable. This is true when a constraint qualification, e.g., Slater's condition, holds.
Let $\lambda ^*$ denote the optimal dual variable.
In the sequel, we invoke Algorithm \ref{DSA_2} to solve the dual problem in \eqref{dual_Lagrangian}. Specifically, we choose
$
d(\lambda)=\frac{1}{2}\lVert \lambda \rVert_2^2
$
as the prox-function of the feasible set. Detailed steps at each round for agents are summarized in Algorithm \ref{consensus-based_dual_decomposition} and explained as follows.

Each agent in Algorithm \ref{consensus-based_dual_decomposition} initializes the algorithm by setting $s_{i,0}=\triangledown \psi_i(\lambda_{i,0})$. This calls for a local maximization step, i.e., $x_i(\lambda_{i,0}) \in \arg \max_{x_i\in\mathcal{X}_i} \big\{-f_i(x_i)-\langle \lambda_{i,0}, h_i(x_i)\rangle \big\}$, to calculate $\triangledown \psi_i(\lambda_{i,0})$ according to Danskin's Theorem. For the initialization of $\lambda_{i,0}$, a reasonable choice is to set $\lambda_{i,0}=0_m$ such that $x_i(\lambda_{i,0}) \in \arg \max_{x_i\in\mathcal{X}_i} -f_i(x_i) $, which is the solution to \eqref{coupled_constraints} if the global constraint is removed. At each round, each agent makes use of ${\rm DSA_2}$ to minimize \eqref{dual_Lagrangian} by performing steps $6$, $7$ and $10$, which correspond to the step $5$ in Algorithm \ref{DSA_2}. To derive $\triangledown \psi_i(\lambda_{i,t}) $, agents should further conduct step $8$. Step $9$ can be seen as a primal recovery step, a common feature that can be found in the literature addressing dual decomposition \cite{simonettoprimalrecovery,falsonedualdecomposition}. This issue is mainly due to the fact that the dual objective function at the optimum is generally nonsmooth, i.e., optimal dual variable does not necessarily lead to an optimal primal solution \cite{nesterovdualsubgradient}. 
\begin{algorithm}
	\begin{algorithmic}[1]
		\caption{${\rm DSA_2}$-based dual decomposition}
		\label{consensus-based_dual_decomposition}
		\STATE  Set $t=0, s_{i,0}=\triangledown \psi_i(\lambda_{i,0})$, choose a non-decreasing sequence of positive parameters $\{\gamma_t\}_{t\geq 0}$.		
		\WHILE  { { Convergence is not reached}}
		\FOR{Each agent $i\in\mathcal{V}$ (in parallel)}  
		\STATE  Receive $s_{j,t}, \forall j\in\mathcal{N}_{i}$;
		\STATE Conduct
		
		\STATE$\hat{\lambda}_{i,t+1}=\arg \min_{\lambda\geq 0}\big\{\big\langle \sum_{k=0}^{t}s_{i,k}, \lambda \big\rangle +\gamma_td(\lambda)\big \}$
		
		\STATE$	\lambda_{i,t+1} = \frac{t+1}{t+2}{\lambda}_{i,t}+\frac{1}{t+2}\hat{\lambda}_{i,t+1}=\frac{1}{t+2}\sum_{k=0}^{t+1}\hat{\lambda}_{i,k}$

		\STATE$
		x_i(\lambda_{i,t+1})=\arg \max_{x_i\in\mathcal{X}_i} -L_i(x_i,\lambda_{i,t+1}) 
		$
		
		\STATE$
		x_{i,t+1}= \frac{t+1}{t+2}{x}_{i,t}+\frac{1}{t+2}x_i(\lambda_{i,t+1})
		$
		
		\STATE$
		s_{i,t+1}=\sum_{j\in\mathcal{N}_i\cup\{i\}}p_{ij}s_{j,t}+\triangledown \psi_i(\lambda_{i,t+1})-\triangledown \psi_i(\lambda_{i,t}) 
		$;
		\STATE Broadcast $s_{i,t+1}$;
		\ENDFOR 
		\STATE Set $t = t+1$.
		\ENDWHILE
	\end{algorithmic}
\end{algorithm}

In the sequel, we establish convergence results for
the dual objective error $\sum_{j=1}^{n}\big(\psi_j(\lambda_{i,t})-\psi_j(\lambda^*)\big)$, the quadratic penalty for the coupled constraint $\big\lVert \big(\sum_{j=1}^{n}h_j(x_{j,t})\big)_+ \big\rVert_2^2$, and the primal objective error $ \sum_{j=1}^{n} \big( f_j(x_{j,t}) -f_j^*\big) $
 for the sequences of test points generated by Algorithm \ref{consensus-based_dual_decomposition} by extending the results derived in previous sections.

\begin{thm}\label{DSA_2_dual_decomposition}
	Suppose Assumption \ref{assumption_weight_matrix} holds true.
	Let the sequences $\{\lambda_{i,t}\}_{t\geq 0}$ and $\{x_{i,t}\}_{t\geq 0}$ be generated by Algorithm \ref{consensus-based_dual_decomposition}. If $\gamma_{t}=\gamma \sqrt{t+1}$ where $\gamma> 0$ and $\lVert \triangledown \psi_i(\lambda_{i,k}) \rVert_2$ is bounded from above, then the dual objective error 
\begin{equation*}
\begin{split}
\sum_{j=1}^{n}\big(\psi_j(\lambda_{i,t})-\psi_j(\lambda^*) \big) 
\leq     \frac{2n\big(\frac{3\sqrt{n}}{1-\sigma_2(P)}+\frac{13}{2}\big)D}{\gamma \sqrt{t+1} }+\frac{\gamma\lVert \lambda^* \rVert_2^2}{2\sqrt{t+1}}
, 
\end{split}
\end{equation*}
the quadratic penalty for the coupled constraint 
\begin{equation*}
	\begin{split}
\big\lVert \big(\sum_{j=1}^{n}h_j(x_{j,t})\big)_+ \big\rVert_2^2\leq \frac{4n\big(\frac{\sqrt{n}}{1-\sigma_2(P)}+\frac{5}{2}\big)D}{t+1}+\frac{2\gamma C}{\sqrt{t+1}},
	\end{split}
	\end{equation*}
	 and the primal objective error \begin{equation*}
	 \begin{split}
	 -\lVert \lambda^*\rVert_2& \sqrt{\frac{4n\big(\frac{\sqrt{n}}{1-\sigma_2(P)}+\frac{5}{2}\big)D}{t+1}+\frac{2\gamma C}{\sqrt{t+1}}   }\\
	 & \leq \sum_{j=1}^{n} \big( f_j(x_{j,t}) -f_j^*\big) 
	  \leq  \frac{2n\big(\frac{\sqrt{n}}{1-\sigma_2(P)}+\frac{5}{2}\big)D}{\gamma \sqrt{t+1} },
	 \end{split}
	 \end{equation*}
	 	where $D= \big(\max_{j\in\{1,\cdots,n\}}\lVert \triangledown \psi_j(\lambda_{j,k}) \rVert_2\big)^2$ and $C=\sum_{j=1}^{n} f_j^*- \min_{\{x_j\in \mathcal{X}_j\}_{j=1}^n}\sum_{j=1}^{n} f_j(x_{j})$ are constants.
\end{thm}

\begin{pf}
In light of \eqref{convexity_1}, we readily have
\begin{equation*}
\begin{split}
&(t+1) \sum_{j=1}^{n}\psi_j(\lambda_{j,t})-\sum_{k=0}^{t}\sum_{j=1}^{n}\psi_j(\lambda_{j,k}) \\
\leq & \sum_{j=1}^{n}\sum_{k=1}^{t}k\langle\triangledown \psi_j(\lambda_{j,k}), \lambda_{j,k}-\lambda_{j,k-1} \rangle .
\end{split}
\end{equation*}
By adding $\sum_{j=1}^{n}\sum_{k=0}^{t}\langle\triangledown \psi_j(\lambda_{j,k}), \lambda_{j,k}-\lambda \rangle, \forall \lambda\geq 0_m $ on both sides, we obtain
\begin{equation} \label{key_dual_lagrangian}
\begin{split}
&(t+1) \sum_{j=1}^{n}\psi_j(\lambda_{j,t})-\sum_{k=0}^{t}\sum_{j=1}^{n}\big( \langle\triangledown \psi_j(\lambda_{j,k}), \lambda -\lambda_{j,k}\rangle\\
&+\psi_j(\lambda_{j,k})  \big)\\ 
\leq&  \sum_{j=1}^{n}\big(\sum_{k=1}^{t}\langle\triangledown \psi_j(\lambda_{j,k}), (k+1)\lambda_{j,k}-k\lambda_{j,k-1}-\lambda \rangle  \\
&+ \langle \triangledown \psi_j(\lambda_{j,0}), \lambda_{j,0}-\lambda \rangle\big) \\
=& \sum_{j=1}^{n}\sum_{k=0}^{t}\langle\triangledown \psi_j(\lambda_{j,k}), \hat{\lambda}_{j,k}-\lambda \rangle,
\end{split}
\end{equation}
where we use the fact that $
(t+1)\lambda_{j,t}=t\lambda_{j,t-1}+\hat{\lambda}_{j,t}
$ to get the last equality.
According to Danskin's Theorem, we have
$
\triangledown \psi_i (\lambda)= -h_i\big(x_i(\lambda)\big) ,
$
where
$
x_i(\lambda) \in \arg \max_{x_i\in\mathcal{X}_i} \big\{-f_i(x_i)-\langle \lambda, h_i(x_i)\rangle \big\}.
$ 
Then by the definition of $\psi_j(\lambda)$, we are able to obtain
\begin{equation*}
\begin{split}
&\sum_{j=1}^{n}\big( \langle\triangledown \psi_j(\lambda_{j,k}), \lambda -\lambda_{j,k}\rangle+\psi_j(\lambda_{j,k})  \big)\\
=& \sum_{j=1}^{n} \Big( \big\langle -h_j\big(x(\lambda_{j,k})\big), \lambda-\lambda_{j,k} \big\rangle -f_j\big(x(\lambda_{j,k})\big)\\
&-\big\langle \lambda_{j,k},h_j\big(x(\lambda_{j,k})\big) \big\rangle\Big) \\
=& \sum_{j=1}^{n}\Big(-\big\langle h_j\big(x(\lambda_{j,k})\big),\lambda\big\rangle-f_j\big(x(\lambda_{j,k})\big) \Big).
\end{split}
\end{equation*}
Plugging the preceding relation into \eqref{key_dual_lagrangian} leads to
\begin{equation}\label{danskin_theorem}
\begin{split}
&\sum_{k=0}^{t}\sum_{j=1}^{n}\Big(f_j\big(x_j(\lambda_{j,k})\big)+\big\langle h_j\big(x_j(\lambda_{j,k})\big),\lambda\big\rangle\Big)\\
&+(t+1) \sum_{j=1}^{n}\psi_j(\lambda_{j,t})
\leq  \sum_{j=1}^{n}\sum_{k=0}^{t}\langle\triangledown \psi_j(\lambda_{j,k}), \hat{\lambda}_{j,k}-\lambda \rangle.
\end{split}
\end{equation}
Recall by definition that $x_{j,t}=\frac{1}{t+1}\sum_{k=0}^{t}x_j(\lambda_{j,k})$. Then by convexity of $f_j$ and $h_j$, we obtain from \eqref{danskin_theorem} that
\begin{equation*}
\begin{split}
&(t+1) \sum_{j=1}^{n} \big( f_j(x_{j,t})+ \langle h_j(x_{j,t}), \lambda \rangle+\psi_j(\lambda_{j,t})\big)\\
\leq&  \sum_{j=1}^{n}\sum_{k=0}^{t}\langle\triangledown \psi_j(\lambda_{j,k}), \hat{\lambda}_{j,k}-\lambda \rangle. 
\end{split}
\end{equation*}
We now follow the same line as in Lemma \ref{basic_convergence_thm} and Lemma \ref{disagreement} to bound the right-hand side of the preceding inequality. Consider
\begin{equation*}
\begin{split}
&(t+1) \sum_{j=1}^{n} \big( f_j(x_{j,t})+\langle h_j(x_{j,t}),\lambda \rangle-(-\psi_j(\lambda_{j,t}))\big)\\
\leq & \sum_{k=0}^{t}  \frac{n}{\gamma_{k-1} } (\frac{\sqrt{n}}{1-\sigma_2(P)}+\frac{5}{2})(\max_{j\in\{1,\cdots,n\}}\lVert \triangledown \psi_j(\lambda_{j,k}) \rVert_2)^2\\
&+ \gamma_{t} d(\lambda).
\end{split}
\end{equation*}
Rearranging the terms yields
\begin{equation*}
\begin{split}
& \sum_{j=1}^{n} \Big( f_j(x_{j,t}) -\big(-\psi_j(\lambda_{j,t})\big)\Big)\\
\leq &  \frac{1}{t+1} \sum_{k=0}^{t}  \frac{n\big(\frac{\sqrt{n}}{1-\sigma_2(P)}+\frac{5}{2}\big)D}{\gamma_{k-1} } \\
& + \min_{\lambda\geq 0_m}\{ \frac{\gamma_{t}}{t+1} d(\lambda)- \langle \sum_{j=1}^{n}h_j(x_{j,t}), \lambda \rangle \} \\
= & \frac{1}{t+1} \sum_{k=0}^{t}  \frac{n\big(\frac{\sqrt{n}}{1-\sigma_2(P)}+\frac{5}{2}\big)D}{\gamma_{k-1} } -  \frac{t+1}{2\gamma_{t}} \big\lVert \big(\sum_{j=1}^{n}h_j(x_{j,t})\big)_+ \big\rVert_2^2, 
\end{split}
\end{equation*}
or, equivalently,
\begin{equation} \label{feasibility_measure}
\begin{split}
& \sum_{j=1}^{n} \Big( f_j(x_{j,t}) -\big(-\psi_j(\lambda_{j,t})\big)\Big) +  \frac{t+1}{2\gamma_{t}} \big\lVert \big(\sum_{j=1}^{n}h_j(x_{j,t})\big)_+ \big\rVert_2^2 \\
&\leq  \frac{1}{t+1} \sum_{k=0}^{t}  \frac{n\big(\frac{\sqrt{n}}{1-\sigma_2(P)}+\frac{5}{2}\big)D}{\gamma_{k-1} }\leq    \frac{2n\big(\frac{\sqrt{n}}{1-\sigma_2(P)}+\frac{5}{2}\big)D}{\gamma \sqrt{t+1} } .
\end{split}
\end{equation}
Following from the saddle point inequality, we obtain
\begin{equation}\label{saddle_point}
\sum_{j=1}^{n}f_j^*\leq \sum_{j=1}^{n}\big(f_j(x_{j,t})+\langle\lambda^*, h_j(x_{j,t})\rangle\big).
\end{equation}

Adding $\frac{t+1}{2\gamma_{t}} \big\lVert \big(\sum_{j=1}^{n}h_j(x_{j,t})\big)_+ \big\rVert_2^2$ and subtracting $ \sum_{j=1}^{n}  \big(-\psi_j(\lambda_{j,t})\big)$ on both sides yield
\begin{equation*}
\begin{split}
&\sum_{j=1}^{n}\Big(f_j^* - \big(-\psi_j(\lambda_{j,t})\big)-\langle\lambda^*, h_j(x_{j,t})\rangle\Big)\\
&+\frac{t+1}{2\gamma_{t}} \big\lVert \big(\sum_{j=1}^{n}h_j(x_{j,t})\big)_+ \big\rVert_2^2 \\
\leq& \sum_{j=1}^{n}\Big(f_j(x_{j,t})- \big(-\psi_j(\lambda_{j,t})\big)\Big)\\
&+\frac{t+1}{2\gamma_{t}} \big\lVert \big(\sum_{j=1}^{n}h_j(x_{j,t})\big)_+ \big\rVert_2^2\\
\leq&     \frac{2n\big(\frac{\sqrt{n}}{1-\sigma_2(P)}+\frac{5}{2}\big)D}{\gamma \sqrt{t+1} }.
\end{split}
\end{equation*}
Since 
\begin{equation*}
\begin{split}
&-\big\langle \lambda^*,\sum_{j=1}^{n}h_j(x_{j,t})\big\rangle+\frac{t+1}{2\gamma_{t}} \big\lVert \big(\sum_{j=1}^{n}h_j(x_{j,t})\big)_+ \big\rVert_2^2\\
&\geq \min_{z\in\mathbb{R}^m}\big\{ -\langle \lambda^*,z\rangle+\frac{t+1}{2\gamma_{t}} \lVert z_+ \rVert_2^2 \big\} \\
& = -\frac{\gamma_t}{2(t+1)}\lVert \lambda^* \rVert_2^2 = -\frac{\gamma\lVert \lambda^* \rVert_2^2}{2\sqrt{t+1}}
\end{split}
\end{equation*}
and $\sum_{j=1}^{n}f_j^*=\sum_{j=1}^{n}-\psi_j(\lambda^*)$, we have
\begin{equation*}
\begin{split}
\sum_{j=1}^{n}\big(\psi_j(\lambda_{j,t})-\psi_j(\lambda^*) \big) 
\leq     \frac{2n\big(\frac{\sqrt{n}}{1-\sigma_2(P)}+\frac{5}{2}\big)D}{\gamma \sqrt{t+1} }+\frac{\gamma\lVert \lambda^* \rVert_2^2}{2\sqrt{t+1}}.
\end{split}
\end{equation*}
Following the similar reasoning in \eqref{identical_variable} to bound $\sum_{j=1}^{n}\big(\psi_j(\lambda_{j,t})-\psi_j(\lambda_{i,t}) \big) $, we further have
\begin{equation*}
\begin{split}
\sum_{j=1}^{n}\big(\psi_j(\lambda_{i,t})-\psi_j(\lambda^*) \big) 
\leq     \frac{2n\big(\frac{3\sqrt{n}}{1-\sigma_2(P)}+\frac{13}{2}\big)D}{\gamma \sqrt{t+1} }+\frac{\gamma\lVert \lambda^* \rVert_2^2}{2\sqrt{t+1}}
. 
\end{split}
\end{equation*}
To establish the upper bound on the violation of the coupled constraint, we consider 
$\forall \lambda\geq 0_m$,
\begin{equation}\label{weak_duality}
f_j^*\geq L_j(x_j^*,\lambda) \geq \min_{x_j\in\mathcal{X}_j} L_j(x_j,\lambda)=-\psi_j(\lambda),
\end{equation}
and therefore
\begin{equation*}
\begin{split}
&  \frac{t+1}{2\gamma_{t}} \big\lVert \big(\sum_{j=1}^{n}h_j(x_{j,t})\big)_+ \big\rVert_2^2\\
& \leq \frac{\sum_{k=0}^{t}  \frac{n\big(\frac{\sqrt{n}}{1-\sigma_2(P)}+\frac{5}{2}\big)D}{\gamma_{k-1} } }{t+1} + \sum_{j=1}^{n} f_j^*- \min_{\{x_j\in \mathcal{X}_j\}_{j=1}^n}\sum_{j=1}^{n} f_j(x_{j}),
\end{split}
\end{equation*}
where \eqref{feasibility_measure} is used. This can be equivalently written as
\begin{equation*}
\begin{split}
& \big \lVert \big(\sum_{j=1}^{n}h_j(x_{j,t})\big)_+ \big\rVert_2^2\\
&\leq \frac{2\gamma_t}{(t+1)^2} \Big(\sum_{k=0}^{t}  \frac{n\big(\frac{\sqrt{n}}{1-\sigma_2(P)}+\frac{5}{2}\big)D}{\gamma_{k-1} } \Big)+\frac{2\gamma_{t}C}{t+1}\\
& \leq \frac{4n\big(\frac{\sqrt{n}}{1-\sigma_2(P)}+\frac{5}{2}\big)D}{t+1}+\frac{2\gamma C}{\sqrt{t+1}}.
\end{split}
\end{equation*}

By Eqs. \eqref{feasibility_measure} and \eqref{weak_duality}, we readily have
\begin{equation*}
\begin{split}
 \sum_{j=1}^{n} \big( f_j(x_{j,t}) -f_j^*\big)\leq    \frac{2n\big(\frac{\sqrt{n}}{1-\sigma_2(P)}+\frac{5}{2}\big)D}{\gamma \sqrt{t+1} }  .
\end{split}
\end{equation*}
Again, by the saddle point inequality \eqref{saddle_point},
the fact
\begin{equation*}
\big(\sum_{j=1}^{n}h_j(x_{j,t})\big)_+\geq \sum_{j=1}^{n}h_j(x_{j,t}),
\end{equation*} 
and $\lambda^*\geq 0 $, one obtains
\begin{equation*}
\begin{split}
&\sum_{j=1}^{n} \big( f_j(x_{j,t}) -f_j^*\big)\geq -\lVert \lambda^*\rVert_2  \big\lVert \big(\sum_{j=1}^{n}h_j(x_{j,t})\big)_+\big \lVert_2
\\
& \geq -\lVert \lambda^*\rVert_2 \sqrt{\frac{4n\big(\frac{\sqrt{n}}{1-\sigma_2(P)}+\frac{5}{2}\big)D}{t+1}+\frac{2\gamma C}{\sqrt{t+1}}   }.
\end{split}
\end{equation*}

This completes the proof.

\end{pf}

%
\begin{remark}
	 By Danskin's Theorem, we have $\triangledown \psi_i (\lambda)= -h_i\big(x_i(\lambda)\big) $. We know from the problem statement that $h_i(x_i)$ is convex and defined on a compact domain $\mathcal{X}_i$. Then it may be not restrictive to assume in Theorem \ref{DSA_2_dual_decomposition} that the subgradient of the dual objective is bounded. This assumption has the same role with the Lipschitz continuity of the objective function used in Lemma \ref{basic_convergence_thm}.
\end{remark}

\section{Simulations}
In this section, we verify our theoretical findings by applying them to a setting where a networked multi-agent system of $n=50$ agents solves a distributed optimization problem with coupled nonlinear convex constraints, i.e.,
\begin{equation*}
\begin{split}
&\min_{x_i\in[0,1]} \sum_{i=1}^{50} c_ix_i \\
&{\rm subject \,\ to} \quad \sum_{i=1}^{50}-d_i\log(1+x_i)\leq -b.
\end{split}
\end{equation*}
Networked-structured optimization problems of this form arise in, for example, plug-in electric vehicles charging  and quality service of wireless networks, and have been used for numerical studies also in \cite{mateosdistributed}.

To solve this problem by ${\rm DSA_2}$, we should first derive its Lagrangian
$
L(\{x_i\}_{i=1}^n,\lambda)
=\sum_{i=1}^{n}L_i({x}_i,\lambda)=\sum_{i=1}^{n}\Big(c_ix_i+\big\langle \lambda, \frac{b}{50}-d_i\log(1+x_i) \big\rangle \Big)
$
and consider the corresponding dual problem
$
\min_{\lambda\geq0} \sum_{i=1}^{n}\psi_i(\lambda),
$  
where $\psi_i(\lambda)= \max_{x_i\in[0,1]}\big\{-c_ix_i-\big\langle \lambda, \frac{b}{50}-d_i\log(1+x_i) \big\rangle\big\}$,
as explained in Section \ref{section5}. 
By Danskin's Theorem, we have
$
\triangledown \psi_i (\lambda)= - \frac{b}{50}+d_i\log\big(1+x_i(\lambda)\big) ,
$
where
$
x_i(\lambda) \in \arg \max_{x_i\in[0,1]} \Big\{-c_ix_i-\big\langle \lambda, \frac{b}{50}-d_i\log(1+x_i) \big\rangle \Big\}
$.

In this simulation, the parameters $c_i$ and $d_i$ for each agent $i\in\{1,\cdots,50\}$ are randomly chosen from a uniform distribution, and $b$ is set as $5$. We use the solver \emph{fmincon} with the interior point algorithm in the Optimization Toolbox to calculate the optimal solution and use it as a reference. The communication topology among agents is characterized by a fixed connected small world graph \cite{wattscollective}, and the weighting matrix $P$ is selected as the Metropolis constant edge weight matrix \cite{xiaodistributed}. The prox-function is chosen as $
d(\lambda)=\frac{1}{2}\lVert \lambda \rVert_2^2.
$
The increasing sequence $\gamma_{t}$ is determined as ${0.2}{\sqrt{t+1}}$.
The initial guesses $\lambda_{i,0}$ for all agents are set as $0$. In the simulation, the bound for the subgradient of the dual objective $\max_{j\in\{1,\cdots,n\}}\lVert \triangledown \psi_j(\lambda_{j,k}) \rVert_2$ is identified as $0.55$ and $\sigma_2(P)$ is calculated as $0.9788$. 
$C$ is estimated as $27.2067$. The optimal dual variable returned by \emph{fmincon} is $0.6419$. Therefore, the theoretical bounds for the primal objective error $\big\lvert \sum_{j=1}^{n} \big( f_j(x_{j,t}) -f_j^*\big) \big\rvert$ and the quadratic penalty for the coupled constraint $\big\lVert \big(\sum_{j=1}^{n}h_j(x_{j,t})\big)_+ \big\rVert_2^2$ are $\max\{ \frac{5.0849\times10^4}{\sqrt{t+1}}, \frac{91.5467}{\sqrt{t+1}}+\frac{6.9856}{(t+1)^\frac{1}{4}}\}$ and $\frac{2.0340\times10^4}{t+1}+\frac{10.8827}{\sqrt{t+1}}$, respectively.

For comparison, we also simulate the consensus-based dual decomposition strategies recently reported in \cite{falsonedualdecomposition, mateosdistributed,simonettoprimalrecovery} in the same network environment. For \cite{falsonedualdecomposition}, according to the sufficient conditions for ensuring convergence developed therein, the stepsize is chosen as $\frac{10}{t+1}$. 
For \cite{simonettoprimalrecovery}, we use a constant stepsize 0.05.
For \cite{mateosdistributed}, we derive the critical feasible consensus stepsize according to Proposition $4$ in \cite{mateosdistributed} as $\sigma = 0.1103$, and using the Slater vector $(\mathbf{1}, 0)$ to get the bound on the optimal dual set as $D=3.3130$. The stepsize is chosen by following the Doubling Trick scheme developed in \cite{mateosdistributed}. The initial guesses $\lambda_{i,0}$ for these two strategies are also set as $0$. 

\begin{figure}[!htb]
	\centering%
	\includegraphics[width=3.5in]{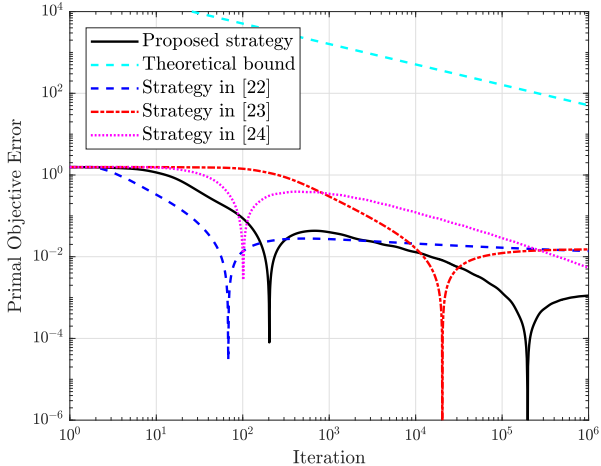}
	\caption{Trajectories of the primal objective error $\big\lvert \sum_{i=1}^{50}c_ix_{i,t} -\sum_{i=1}^{50}c_ix^*_{i}\big\rvert$.} \label{cost_error}
\end{figure}

\begin{figure}[!htb]
	\centering%
	\includegraphics[width=3.5in]{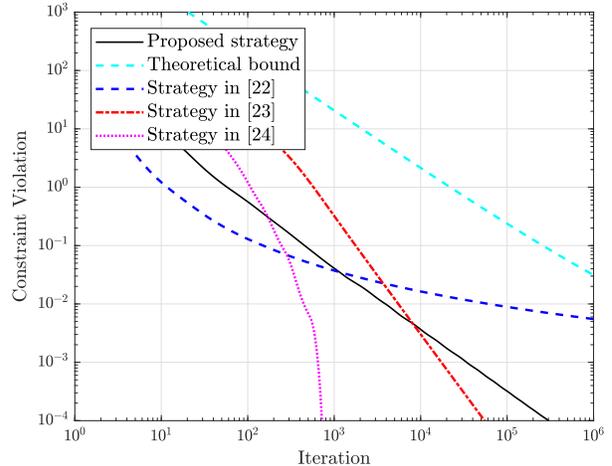}
	\caption{Trajectories of the quadratic penalty for the coupled constraint $\big\lVert\big(b-\sum_{i=1}^{50}d_i\log(1+x_{i,t})\big)_+\big\lVert_2^2$.} \label{cons_violation}
\end{figure}

The simulation results are illustrated in Figs. \ref{cost_error} and \ref{cons_violation}. In particular,
 Figs. \ref{cost_error} and \ref{cons_violation} depict the primal objective error and the quadratic penalty for the coupled constraint, respectively.
 It is worth mentioning that, for methods in \cite{falsonedualdecomposition, mateosdistributed,simonettoprimalrecovery}, the performances are evaluated over the running sequences of the primal variables to be compatible with the theoretical results therein. We can see from them that the algorithm in \cite{falsonedualdecomposition} suffers from the slowest convergence rate among the three, and the proposed algorithm shows a slightly quicker convergence rate than that in \cite{mateosdistributed}.
This may be because that the stepsizes for \cite{mateosdistributed} (Doubling Trick scheme) and the proposed one ($\frac{1}{\gamma_t}=\frac{5}{\sqrt{t+1}}$) are of order $\frac{1}{\sqrt{t}}$ while the stepsize for \cite{falsonedualdecomposition} is chosen to be of order $\frac{1}{t}$ to fulfill the sufficient conditions developed therein.
Due to using a constant stepsize, the method in \cite{simonettoprimalrecovery} does not provide exact convergence.
In Fig. \ref{cost_error}, we draw the absolute value of $ \sum_{i=1}^{50}c_ix_{i,t} -\sum_{i=1}^{50}c_ix^*_{i}$. The identity $\sum_{i=1}^{50}c_ix_{i,t} -\sum_{i=1}^{50}c_ix^*_{i}$ can be negative and positive due to possible violation of the primal coupled constraint. When it jumps from negative to positive, the trajectory of $ \big\lvert \sum_{i=1}^{50}c_ix_{i,t} -\sum_{i=1}^{50}c_ix^*_{i}\big\rvert$ presents a peak.
This phenomenon is typically observed in dual Lagrangian problems.
We also note that the trajectories for the proposed algorithm are within the theoretically developed upper bounds.


\section{Conclusion}
In this paper, we have proposed a distributed subgradient method with double averaging, termed as ${\rm DSA_2}$, for convex constrained optimization problems with aggregate nonsmooth objective functions in a multi-agent setting. In particular, the local iteration rule in ${\rm DSA_2}$ works in the dual space, that is, it involves at each round minimizing a local approximated dual model of the overall objective where the estimated first-order information of the objective is supplied by a dynamic average consensus scheme. A non-ergodic convergence rate of $O(\frac{1}{\sqrt{t}})$ in terms of objective function error  has been established for ${\rm DSA_2}$. Furthermore, we have developed a ${\rm DSA_2}$-based dual decomposition strategy for solving distributed optimization problems with coupled functional constraints. 
This is made possible by dualizing the coupled constraint via Lagrangian relaxation, thus allowing us to alternatively solve the dual Lagrangian problem where the coupling takes place in cost functions via ${\rm DSA_2}$.
The $O(\frac{1}{\sqrt{t}})$ convergence rate has been theoretically validated for both the dual objective error and the quadratic penalty for the coupled constraint. 
Several simulations and comparisons have been performed to verify the advantages of the proposed methods.




\appendix
\end{document}